# A Multidimensional Artistic Approach to Enhance Understanding of Julia Sets through Computer Programming

E. A. Navas-López

Department of Mathematics, University of El Salvador, El Salvador

## Abstract

*This article proposes an artistic approach to increase and enrich the understanding of Julia Sets. This approach includes the mathematical, the playful, the artistic and the computational dimensions. It is argued that these four dimensions are not disjointed or dissociated despite general rejection by traditional academic communities and art critics communities. Also, some significant collections of Computational Art or Computer-Generated Mathematical Art are mentioned. Four artistic creations based on Julia Sets are presented as examples using the CFDG language. Finally, an informal application of the approach was carried out and artistic production of students are presented and discussed.*

## 1. Introduction

There is a relatively recent tendency to incorporate art into the teaching of mathematics that criticizes mathematics education in the sense of overcoming the perspective of traditional teaching [1]. However, over the years, there has been rejection of incorporating art into the mathematical world [2], and there has also been rejection of incorporating computer-generated art into the art world [3].

But, as Ferreira and Lessa [1], the author believes in the potential of Art as an element of struggle and resistance, in the sense of breaking with the reified reality, pointing to horizons of transformation. So, as called by [1], this paper proposes some ideas for the dissemination of educational practices that mobilize the Art and Mathematics interface, also considering the discussions provoked.

Section 1.1 presents a general background reference on the integration of mathematic, artistic, playful and computational dimensions, particularly on Julia sets, computer-generated art, and education. Then, in section 1.2, the mathematical definitions necessary to develop this proposal are presented. Section 1.3 briefly explains how to generate graphical representations of Julia sets using the CFDG language.

Section 2 presents four examples of computer-generated artworks based on Julia sets to use as a starting point for the proposal. Sections 3 and 4 explain the applied methodology and some of the empirical results obtained. Finally, section 5 presents the discussion of the results.

 

## 1.1. Background

As Bergweiler [4] explains, the study of iteration theory is fundamental in mathematics, and its classic problem is the study of the iterative behaviour of a family of functions that depend on a parameter. In this sense appears the study of what we now know as Julia Sets in the early twentieth century. Of course, at the time there were no computers and the study of these sets was very difficult. However, Julia sets play a critical role in the understanding of the dynamics of families of mappings [5].

In recent decades, as Hitt [6, p. 214] points out: "Technological advancement has significantly influenced the development of theoretical notions that were previously taken into account but were not considered crucial in terms of explaining the learning of mathematical concepts. These theoretical aspects are the basis for understanding the study of the different representations of mathematical objects and their role in the construction of concepts. Now, with technology, it is important to study the different representations of mathematical objects in environments very different from those that were followed in the past".

The study of many areas of mathematics and mathematics education have been modified with the popularization of different technologies such as personal computers, including Fractal Geometry. So, building multiple computer-generated images to form a richer mental image ("mental image" in the Vinner's [7] sense) of fractals like Julia Sets is quite affordable for the students of our day. One way to implement this construction to achieve a better understanding of Julia Sets (and other types of fractals) is by artistic means, motivating the learner in a purely playful process to build (in the sense of creating) images not only aesthetically pleasing but also endowed with some meaning through computer programs.

This pathway is not widely used in our time in mathematics in general, because although mathematics and art have been very close since the first manifestations of rationality of the human species, unfortunately, we have seen that these two areas of knowledge have distanced in school programs [8]. In his book, D'Amore [9] makes a brilliant exposition of the presence of mathematics at the dawn of humanity.

The playful dimension is more common. Bishop [10, p. 21] discusses the role of games in mathematics education and notes: "Educators in mathematics have discovered through their experience, and they have supported with theoretical research, that playing can be an integral part of learning. This has made the act of playing and the idea of gaming a much more widespread teaching and learning activity than it had been before". For further reading, in the literature review made by González Peralta, et al. [11], one can find possible research lines about gaming in mathematics education.

There is also the computational dimension, which is becoming more common. For example, Hoffmann [12] presents a sixth-grade primary experience using a Monte Carlo simulation for introducing the concept of area of a unit circle (which is the approximation of the number $\pi$). There is also the Experience of DeJarnette [13] in which students use the Scratch programming environment to help themselves understanding of how distances travelled by certain objects are functions of time.

So, we have the seemingly strange conjunction of the four dimensions: mathematic, artistic, playful and computational. This conjunction is not directly welcomed by mathematicians, artists and art critics, or traditional computer scientists. In fact, pioneers in Computational Art faced rejection by the mathematical community, as Mumford, et al. [2, p. viii] emphasizes: "What to do with the pictures? [...] they were unpublishable in the standard way. There were no theorems, only very suggestive pictures. They furnished convincing evidence for many conjectures and lures to further exploration, but theorems were the coin of the realm and the conventions of that day dictated that journals only publish theorems".

A similar marginalization (and/or misunderstanding) occurs to Rani and Kumar [14] whose article was published under Series D journal, "Research in Mathematical Education" and not under Series B journal, "The Pure and Applied Mathematics" of the same Organization (http://www.ksme.info/eng/), despite their article is not about education at all, but about Superior Mandelbrot Set. This paper is clearly about Pure Mathematics even though there are no new theorems. It has "only new pictures and conjectures".

There are also problems with artistic acceptance according to Franke [3, p. 186]: "[Images] were considered [only] drawings from the plotter, the main problem was the uncertainty of the experts, the art historians and the critics, and above all the gallery owners. The problem was that the computer can produce an arbitrary number of equally good 'originals', which can be a detriment in the business world of art".

On the apparent incompatibility between science and art, and referencing the Peitgen and Richter's book [15], Eilenberg [16, p. 175] offers a conciliation: "It is rather unusual for natural [physical and mathematical] scientists to endeavour with such tenacity to bring their results and insights to the general public, [...] Instead of giving an abstract presentation in so many dry words, they have chosen pictures with a direct, universal appeal – a combination of mathematics and art!". That is, scientific art (whether mathematical art or computational art) can be used primarily to publicize results to the general public.

On the confusion between Computational Art and Standard Art, Franke [3, p. 187] proposes the following reflection: "The art of every age has used the means of its time to give form to artistic innovation. [...] Why shouldn't the computer, that universal medium of information and communication which has even invaded our private homes [and our lives], be used as a medium and instrument of art?". Moreover, Zaleski Filho [8] tells us: "Thus, true art, which has no end in any of its external realizations, has as its identification a spiritual principle that enlivens and surpasses all of them". So, after all, computational art, mathematical art and scientific art in general, as well as all other types of art, need no more justification than their ability to encourage human beings to rejoice in the artworks themselves.

This is how various collections and producers of Computer Generated Mathematical Art have proliferated, such as the Peitgen and Richter's [15] collection which includes many fractal graphics of complex dynamic systems, such as The Bridges Organization [17] which annually holds an international exhibition and competition for mathematical art (not just computer-generated), such as Aslaksen's [18] collection of university courses about art and mathematics, such as the computer-generated mathematical art collection presented by Navas-López [19] which includes various types of fractals and other types of computer graphing techniques, etc. (explained in their catalogue [20]). There are also other less formal but not less impressive collections, such as Nylander's [21], The Context Free Art [22] community gallery, and Math Munch's [23] blog which included not only mathematical art but many interesting things. There are some collections that are very specific like Ross' [24] mathematical analyses from the Sacks number spiral.

Speaking now of mathematical art in the classroom and in the curriculum in general, the author supports the reasoning of Figueiras, et al. [25, p. 46]: "In the Mathematics that are taught, those that in Obligatory Education it is said that they will serve to acquire what is needed in life, who can deny a place to beauty? Do we intend to let ourselves be carried away only by the dubious pragmatism of a mathematics cut short both in the time available for its teaching and in the potentiality of its values?". Moreover, "The study of fractals is a motivating element in students, due to the implicit aesthetics in their constructions and the suggestive that their designs can be" [25, p. 47].

For Bosque, et al. [26, p. 20]: "[...] Aesthetics transcends the Philosophy of Art finding a place in the Philosophy of Science. In this way, if aesthetics is a component of knowledge and also has a role in scientific and mathematical activity, it makes sense to study its influence on the teaching and learning of mathematics".

The proposed Fractal Geometry Activities in the Redondo and Haro [27, 28] High School Classroom are an excellent source of ideas for planning outreach activities to different types of fractals, including Julia Sets (see [28, p.17]). However, in this vast and general overview, the artistic approach and colours are scarce. The presentation of the Mandelbrot and Julia Sets by Varona [29] explains technical details about how to graph them in Mathematica(R), which is proprietary software, and proposes the inclusion of colour palettes to enhance the images although it can only display them in grey-scale due to the type of publication. This is an example that although there is acceptance of the subject, not all journals, editors or publishers are interested or prepared to accept "mathematics with colours".

So, this paper then aims to expand the part of Redondo and Haro's [27, 28] proposal on the Julia Sets topic, using computer programming as in the Varona's [29] paper, but with the following differences: (a) emphasizing the artistic approach and not the mathematical one, (b) not reducing to the aesthetic dimension but incorporating the communicative and didactic dimensions, (c) adding colour, and (d) using free software, in this case the ContextFree software (www.contextfreeart.org).

### 1.2. Julia Sets

We can take the definition of Julia Set from [4, p. 153]: Let $f : \mathbb{C} \to \hat{\mathbb{C}}$ be a meromorphic function, where $\mathbb{C}$ is de complex plane and $\hat{\mathbb{C}} = \mathbb{C} \cup \{\infty\}$. [...] we shall always assume that $f$ is neither constant nor a linear transformation. Denote by $f^n$ the $n$th iterate of $f$, that is, $f^0(z) = z$ and $f^n(z) = f(f^{n-1}(z))$ for $n \geq 1$. The basic objects studied in iteration theory are the Fatou set $F = F(f)$ and the Julia set $J = J(f)$ of a meromorphic function $f$. Roughly speaking, the Fatou set is the set where the iterative behaviour is relatively tame in the sense that points close to each other behave similarly, while the Julia set is the set where chaotic phenomena take place. The formal definitions are:

$F = \{z \in \mathbb{C} : \{f^n : n \in \mathbb{N}\}$ is defined and normal in some neighbourhood of $z\}$
and $J = \hat{\mathbb{C}} - F$.

However, here will be chosen the most simplified version of [30, p. 263]: "Julia set of a function $f_c(z)$ with seed $c \in \mathbb{C}$, denoted by $J_c(f)$, is the set composed by all $z \in \mathbb{C}$, such that the sequence $\{z_k\}_{k=0}^{\infty}$ is bounded, where $z_0 = z$ and $z_{n+1} = f_c(z_n)$ [for $n \geq 0$]. Where $J_c(f)$, in its simplest form, uses $f_c(z) = z^2 + c$". Nevertheless, many other different functions can be used that provide interesting results such as those presented by Entwistle [31], Garijo, et al. [32], Liu, et al. [33], Peitgen and Richter [15], Pickover and Khorasani [34], Rochon [35] and Rani and Kumar [14].

The criterion used to determine whether the sequences diverge is if $z_k$, for some $k \geq 0$, has a modulus greater than 2, that is $|z_k| \geq 2$. Since it cannot be evaluated to infinity, a bound is used: $N$. If the sequence does not "diverge" before reaching the $N$th term (when $k \leq N$), it is considered not divergent, that is, it is bounded. The larger $N$, the greater the precision of the set [30, p. 263].

### 1.3. Implementation of Julia sets pictures on CFDG language

The CFDG language, version 3, of Context Free (https://www.contextfreeart.org/) software and adapted computer graphics techniques from [36] and some from [37] (from their respective chapters on fractal graphing) are used for implementation here. No extra libraries are needed.

The CFDG language is not a programming language properly. It is actually a language in which can be defined a particular type of context-free grammars whose terminal symbols are primitive figures: squares, circles and triangles. Different related transformations (displacement, scaling, rotation, etc.) can be applied to these figures. However, CFDG language supports Functional Programming when it is needed, especially for numerical algorithms. For further reading, visit the Context Free Art documentation page [38].

The basic geometric object to use in the examples is the square, that will be each pixel of images. This is constructed using the following primitive in CFDG language:

```
SQUARE [
    x <x_offset>
    y <y_offset>
    size <width><height>   # Can be abbreviated as s
    hue <tone>             # Can be abbreviated as h
    saturation <saturation> # Can be abbreviated as sat
    brightness <bright>    # Can be abbreviated as b
]
```

Where `<x_offset>` and `<y_offset>` indicate the displacement from the origin, `<width>` and `<height>` determine the size of the figure, `<tone>` is an angle between 0° and 360° indicating the colour of the figure according to the HSV colour model, and `<saturation>` and `<bright>` indicate the corresponding. For more information about HSV colour model, see [39].

The function required to determine the convergence of a point (`z_r`, `z_i`) in the complex plane in CFDG language is in listing 1 (lines 3–8). Where `MAXSTEPS` is the value of $N$ in the criterion described in section 1.2. The initial call must be in the form `steps(0, z_r, z_i, c_r, c_i)`, where (`c_r`, `c_i`) is the seed value. This call returns the number of iterations executed from which the sequence diverges, or returns `MAXSTEPS` if the sequence has not yet diverged at the $N$th iteration.

Listing 1: Basic source code for Julia set in CFDG language

```
1. startshape julia(-0.381966, 0.618034)
2.
3. MAXSTEPS = 40
4. steps(numSteps, z_r,z_i, c_r,c_i) =
5.     if((numSteps < MAXSTEPS) && (z_r*z_r+z_i*z_i<4),
6.         steps(numSteps+1,
7.             z_r*z_r - z_i*z_i + c_r, 2*z_r*z_i + c_i, c_r, c_i),
8.         numSteps)
9.
10. LIMIT = 1000 # Image resolution
11.
12. # Borders of the complex plane to show:
```

```
13. LIMLEFT = -1.4
14. LIMRIGHT = 1.4
15. LIMTOP = 1.4
16. LIMBOT = -1.4
17.
18. # Width and height of the squares that will discretize the image:
19. SIZEX = (LIMRIGHT-LIMLEFT)/(LIMIT-1)
20. SIZEY = (LIMTOP-LIMBOT)/(LIMIT-1)
21.
22. shape julia(c_r,c_i) {
23.     loop i = (LIMIT) [] {
24.         z_i = (LIMTOP-LIMBOT)*i/(LIMIT-1) + LIMBOT # y
25.         loop j = LIMIT [] {
26.             z_r = (LIMRIGHT-LIMLEFT)*j/(LIMIT-1) + LIMLEFT # x
27.
28.             numSteps = steps(0, z_r, z_i, c_r, c_i)
29.             if (numSteps==MAXSTEPS){
30.                 # Black
31.                 SQUARE[x z_r y z_i size SIZEX SIZEY b 0]
32.             } else {
33.                 # Gray
34.                 SQUARE[x z_r y z_i size SIZEX SIZEY b 0.9]
35.             }
36.         }
37.     }
38. }
```

Then a Julia set figure must be implemented, determined by the seed, with the definitions in lines 10–38 from listing 1. The `startshape` directive is used to indicate which is the generating/starting shape (see listing 1, line 1).

Result of execution of listing 1 source code with `MAXSTEPS` values equal to 40, 60, 80 and 100 is presented in figure 1.

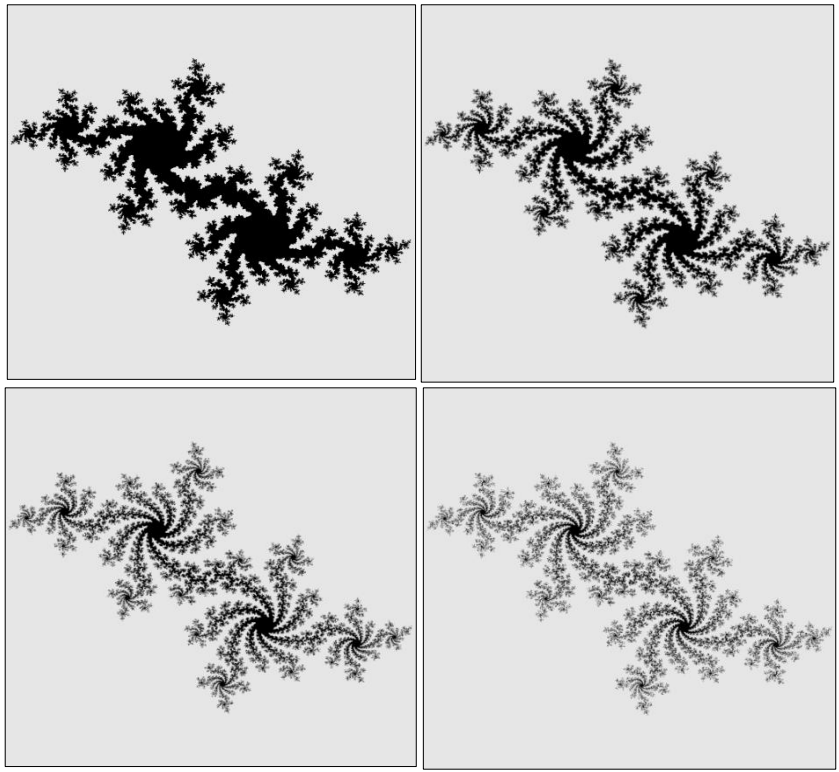

Figure 1. Basic Julia set image generated by CFDG code, with values of $N$=40, 60, 80 and 100.

## 2. ARTWORKS RAISED AS EXAMPLES

It is a recommended collection of examples as a proposal for experiment with the parameters like seed number, viewport (drawn interval), colours, bright, saturation and hue formulae, etc.

### 2.1. Frozen Fjords

This artwork (Figure 2) shows an aerial view of snow-capped fjords, its thin dark sand shores and the deeply blue sea. The motivation is that fjords have a natural fractal shape.

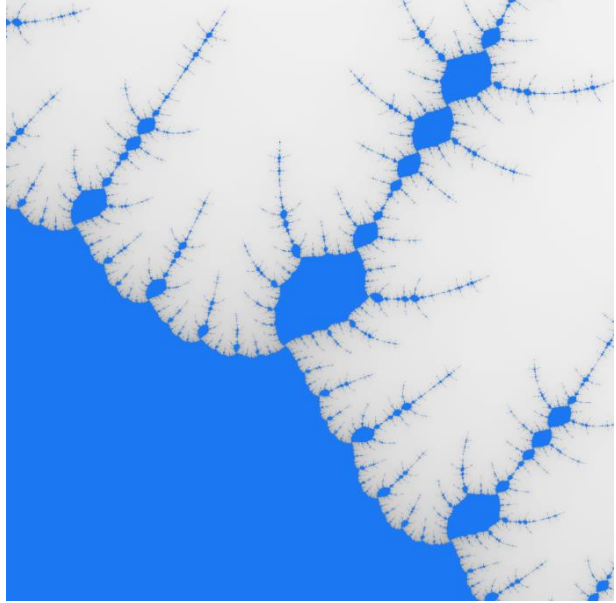

Figure 2. Frozen Fjords

From a technical point of view, it is a view of a Julia set in the range of [0.01, 0.09]×[0.02$i$, 0.1$i$], with see $c = -1.384286+0.004286i$. The colouration of this artwork has a constant hue as well as saturation, but the brightness is variable depending on the number of steps in which it is determined that the point belongs or does not belong to the set. See Listing 2.

Listing 2: CFDG source code for Frozen Fjords

```
1.  # This file is Free Software released under the GNU GPLv3 license or
2.  # its latest version:
3.  # http://www.gnu.org/licenses/gpl.html
4.
5.  # To generate the image run the following line:
6.  # $ cfdg -b 0 -s 1000 fjords.cfdg fjords.png
7.
8.  startshape fjords(-1.384286,0.004286)
9.
10. LIMIT = 1000 # Image resolution
11. MAXSTEPS = 300
12.
13. # Borders of the complex plane to show:
14. LIMLEFT = 0.01
15. LIMRIGHT = 0.09
16. LIMTOP = 0.10
17. LIMBOT = 0.02
18.
19. # Width and height of the boxes that will discretize the image:
20. SIZEX = (LIMRIGHT-LIMLEFT)/(LIMIT-1)
```

```
21. SIZEY = (LIMTOP-LIMBOT)/(LIMIT-1)
22.
23. steps(numSteps,z_r,z_i,c_r,c_i) =
24.       if((numSteps < MAXSTEPS) && (z_r*z_r+z_i*z_i<4),
25.           steps(numSteps+1,
26.           z_r*z_r - z_i*z_i + c_r, 2*z_r*z_i + c_i, c_r, c_i),
27.           numSteps)
28.
29. shape fjords(c_r,c_i) {
30.     FILL[h 214 sat 0.89 b 0.95] # Blue ocean
31.     loop i = (LIMIT) []  {
32.         z_i = (LIMTOP-LIMBOT)*i/(LIMIT-1) + LIMBOT # y
33.         loop j = LIMIT [] {
34.             z_r = (LIMRIGHT-LIMLEFT)*j/(LIMIT-1) + LIMLEFT # x
35.
36.             numSteps = steps(0, z_r, z_i, c_r, c_i)
37.             if(numSteps<MAXSTEPS){
38.                 SQUARE[x z_r
39.                     y z_i
40.                     size SIZEX SIZEY
41.                     h 30 sat 0
42.                     b (1+(1-numSteps)/(MAXSTEPS-1))]
43.             }
44.         }
45.     }
46. }
```

## 2.2. The Wail of the Pripyat Forest

This artwork (Figure 3) shows a sick forest around the city of Pripyat. This ghost town is known for being affected by the worst accident in nuclear power history on April 26, 1986, when the Chernobyl Nuclear Power Plant reactor number 4 was overheated and exploded during a shutdown test. The motivation came after seeing a detailed documentary about the nuclear disaster at the Chernobyl nuclear power plant.

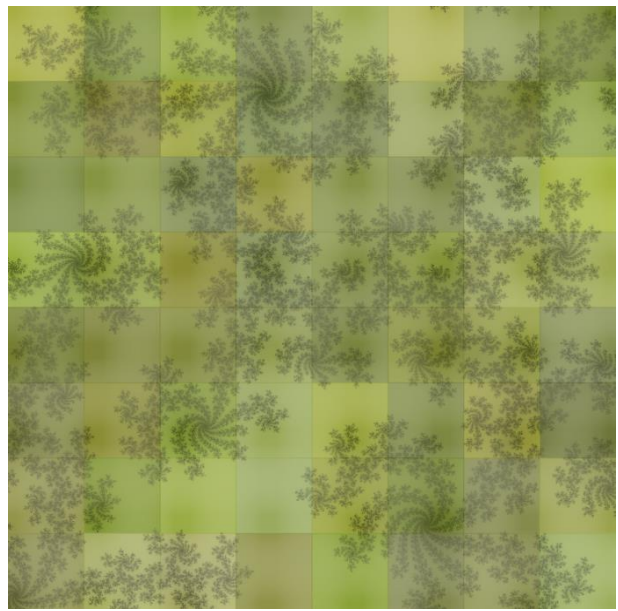

Figure 3. The Wail of the Pripyat Forest

From a technical point of view, it is a view of a Julia set in the range $[-0.052857, 0.188571] \times [-0.105714i, 0.135714i]$ with see $c = -0.381966 + 0.618034i$. The colouration of this artwork has a random variable hue between 60 and 74, a variable random saturation between 0.41 and 0.66, and a variable random maximum brightness between 0.32 and 0.35 for points that do not belong to the set. See Listing 3.

Listing 3: CFDG source code for Wail of the Pripyat Forest

```
1. # This file is Free Software released under the GNU GPLv3 license or
2. # its latest version:
3. # http://www.gnu.org/licenses/gpl.html
4.
5. # To generate the image run the following line:
6. # $ cfdg -b 0 -s 1000 -v PAJBHA forest.cfdg forest.png
7.
8. LIMIT = 1000 # Image resolution
9. MAXSTEPS = 200
10. MINMAXBRIGHT = 0.32
11. MAXMAXBRIGHT = 0.68
12.
13. startshape  forest
14.
15. # Borders of the complex plane to show:
16. LIMLEFT = -0.052857
17. LIMBOT = -0.105714
18. LIMRIGHT = 0.188571
19. LIMTOP = 0.135714
20.
21. # Width and height of the boxes that will discretize the image:
22. SIZEX = (LIMRIGHT-LIMLEFT)/(LIMIT-1)
23. SIZEY = (LIMTOP-LIMBOT)/(LIMIT-1)
24.
25. steps(numSteps,z_r,z_i,c_r,c_i) =
26. if((numSteps < MAXSTEPS) && (z_r*z_r+z_i*z_i<4),
27. steps(numSteps+1,
28. z_r*z_r - z_i*z_i + c_r, 2*z_r*z_i + c_i, c_r, c_i),
29. numSteps)
30.
31. NUMBLOCKS = 8 # The number of columns and rows
32.
33. shape forest{
34.     loop i = NUMBLOCKS []{
35.         loop j = NUMBLOCKS [] {
36.             julia2(
37.                 -0.381966, 0.618034,
38.                 j*(LIMRIGHT-LIMLEFT)/NUMBLOCKS+LIMLEFT,
39.                 i*(LIMTOP-LIMBOT)/NUMBLOCKS+LIMBOT,
40.                 rand(MINMAXBRIGHT,MAXMAXBRIGHT)
41.             ) [h rand(60,74)
42.             sat rand(0.41,0.66)
43.             b rand(0.32,0.35)]
44.         }
45.     }
46. }
47.
48. shape julia2(c_r,c_i, xi,yi, maxBright) {
49.     xf = xi+(LIMRIGHT-LIMLEFT)/NUMBLOCKS
50.     yf = yi+(LIMTOP-LIMBOT)/NUMBLOCKS
51.     loop i = LIMIT/NUMBLOCKS [] {
52.         z_i = (yf-yi)*i/(LIMIT/NUMBLOCKS-1) + yi
53.         loop j = LIMIT/NUMBLOCKS [] {
54.             z_r = (xf-xi)*j/(LIMIT/NUMBLOCKS-1) + xi
55.
56.             numSteps = steps(0, z_r, z_i, c_r, c_i)
57.             SQUARE[x z_r y z_i size SIZEX SIZEY
58.             b ( maxBright+maxBright*(1-numSteps)/(MAXSTEPS-1) )]
59.         }
60.     }
61. }
```

## 2.3. Ragnarök

This artwork (Figure 4) shows a rough, sharp and Nordic scenario illustrating the cataclysm of the Ragnarök, which is the Viking apocalypse. The motivation comes from the story of the violent Viking apocalypse.

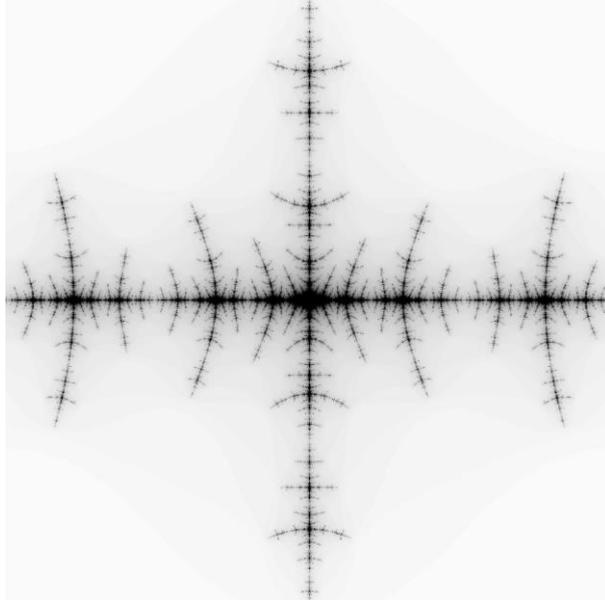

Figure 4. Ragnarök

From a technical point of view, it is a view of a Julia set in the range of $[-0.6, 0.6] \times [-0.6i, 0.6i]$, with seed $c = -1.4+0.0i$. The colouration of this artwork has a higher variable white brightness the faster it is determined that the dots do not belong to the set. The points that do belong to the set are black. See listing 4.

Listing 4: Source code for Ragnarök

```
1.  # This file is Free Software released under the GNU GPLv3 license or
2.  # its latest version:
3.  # http://www.gnu.org/licenses/gpl.html
4.
5.  # To generate the image run the following line:
6.  # $ cfdg -b 0 -s 1000 ragnarok.cfdg ragnarok.png
7.
8.  startshape ragn(-1.4, 0.0)
9.
10. LIMIT = 1000 # Image resolution
11. MAXSTEPS = 100
12.
13. # Borders of the complex plane to show:
14. LIMLEFT = -0.6
15. LIMRIGHT = 0.6
16. LIMTOP = 0.6
17. LIMBOT = -0.6
18.
19. # Width and height of the boxes that will discretize the image:
20. SIZEX = (LIMRIGHT-LIMLEFT)/(LIMIT-1)
21. SIZEY = (LIMTOP-LIMBOT)/(LIMIT-1)
22.
23. steps(numSteps,z_r,z_i,c_r,c_i) =
24.     if((numSteps < MAXSTEPS) && (z_r*z_r+z_i*z_i<4),
25.         steps(numSteps+1,
```

```
26.             z_r*z_r - z_i*z_i + c_r, 2*z_r*z_i + c_i, c_r, c_i),
27.         numSteps)
28.
29. shape ragn(c_r,c_i) {
30.     loop i = LIMIT/2 [] {
31.         z_i = (LIMTOP-LIMBOT)*i/(LIMIT-1) + LIMBOT # y
32.         loop j = LIMIT/2 [] {
33.             z_r = (LIMRIGHT-LIMLEFT)*j/(LIMIT-1) + LIMLEFT # x
34.
35.             numSteps = steps(0, z_r, z_i, c_r, c_i)
36.             bright = (1+(1-numSteps)/(MAXSTEPS-1))
37.             # Symmetry when imag part of seed is zero and viewport is centred:
38.             SQUARE[x z_r    y z_i    size SIZEX SIZEY b bright]
39.             SQUARE[x (-z_r) y z_i    size SIZEX SIZEY b bright]
40.             SQUARE[x z_r    y (-z_i) size SIZEX SIZEY b bright]
41.             SQUARE[x (-z_r) y (-z_i) size SIZEX SIZEY b bright]
42.         }
43.     }
44. }
```

## 2.4. The Battle for Smolensk

Here (Figure 5) what I want to represent is a blood bath on the icy ground of Smolensk at the end of 1941, early winter. The motivation comes from a series of documentaries about World War II, particularly about Operation Barbarossa and how the facts of the German advance on Soviet ground developed in the period 1941-1943.

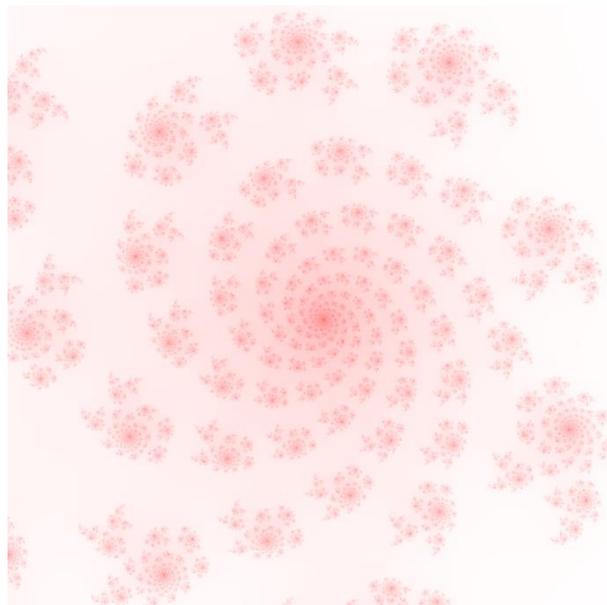

Figure 5. The Battle for Smolensk

From a technical point of view, it is a view of a Julia set in the range of $[-0.21, 0.63] \times [-0.865714i, -0.025714i]$, with center $0.21 - 0.445714i$ and seed $c = 0.39 - 0.252857i$. The colouration of this artwork are red dots and maximum brightness with lower saturation for the points the faster it is determined that they do not belong to the set, and the higher the longer it takes to determine not belonging. The dots that do belong to the set are painted as intense red.

Listing 5: CFDG source code for The Battle for Smolensk

```
1.  # This file is Free Software released under the GNU GPLv3 license or
2.  # its latest version:
3.  #   http://www.gnu.org/licenses/gpl.html
4.
5.  # To generate the image run the following line:
6.  # $ cfdg -b 0 -s 1000 battle.cfdg battle.png
7.
8.  LIMIT = 1000 # Image resolution
9.  MAXSTEPS = 400
10.
11. startshape julia3(0.39, -0.252857)
12.
13. SIDE = 0.84 # side of the viewport square
14. CX = 0.21    # center x
15. CY = -0.445714 # center y
16.
17. # Borders of the complex plane to show:
18. LIMLEFT = CX - SIDE/2
19. LIMBOT = CY - SIDE/2
20. LIMRIGHT = CX + SIDE/2
21. LIMTOP = CY + SIDE/2
22.
23. # Width and height of the boxes that will discretize the image:
24. SIZEX = (LIMRIGHT-LIMLEFT)/(LIMIT-1)
25. SIZEY = (LIMTOP-LIMBOT)/(LIMIT-1)
26.
27. steps(numSteps,z_r,z_i,c_r,c_i) =
28.     if((numSteps < MAXSTEPS) && (z_r*z_r+z_i*z_i<4),
29.         steps(numSteps+1,
30.         z_r*z_r - z_i*z_i + c_r, 2*z_r*z_i + c_i, c_r, c_i),
31.         numSteps)
32.
33. shape julia3(c_r,c_i) {
34.     loop i = LIMIT [] {
35.         z_i = (LIMTOP-LIMBOT)*i/(LIMIT-1) + LIMBOT # y
36.         loop j = LIMIT [] {
37.             z_r = (LIMRIGHT-LIMLEFT)*j/(LIMIT-1) + LIMLEFT # x
38.
39.             numSteps = steps(0, z_r, z_i, c_r, c_i)
40.             SQUARE[x z_r y z_i size SIZEX SIZEY b 1
41.                 sat ((numSteps-1)/(MAXSTEPS-1))]
42. }
43. }
44. }
```

## 3. METHODOLOGY

An informal 8-hour on-line course (4 Saturdays) was implemented with a small group of university students, where the basic elements of the CFDG language were studied. See [38]. In addition, the examples in section 2 were explained. They went back to these examples and made modifications to the source code according to their own creative process or parameter changing experimentation, without a plug-in formula or plug-in solution (in the sense of [40]).

After the 6-hour course (3 Saturdays), they were asked in the fourth session (last one) to make their own artwork based on Julia sets and using CFDG language. The best results are presented in the next section.

## 4. SOME RESULTS

The students experimented a lot with the basic elements and then combined much of them to achieve some interesting patterns. However, only the best results of the final activity are shown here.

### 4.1. Under the shade of leaves

This artwork (figure 6) shows a view from under some leaves. As the author say "I imagined myself lying under a tree looking up at the clear sky."

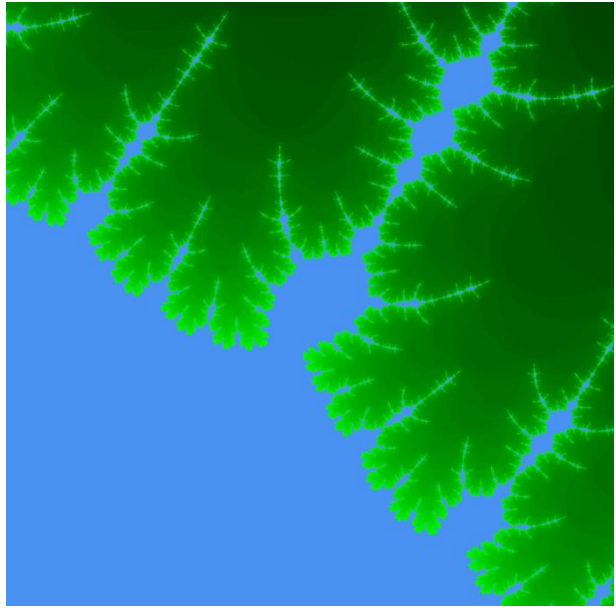

Figure 6. Under the shade of leaves

He based on source code from listing 2, and changed the colours and the brightness calculation formula. He also reduced the maximum number of steps (MAXSTEPS) to make borders more similar to leaf borders. See listing 6.

Listing 6: CFDG source code for Under the shade of leaves

```
1. # Under the shade of leaves.
2. # To generate the image run the following line:
3. # $ cfdg -b 0 -s 1000 leavesshade.cfdg leavesshade.png
4.
5. startshape leaves(-1.384286,0.004286)
6.
7. LIMIT = 1000 # Image resolution
8. MAXSTEPS = 60
9.
10. # Borders of the complex plane to show:
11. LIMLEFT = 0.01
12. LIMRIGHT = 0.09
13. LIMTOP = 0.10
14. LIMBOT = 0.02
15.
16. # Width and height of the boxes that will discretize the image:
17. SIZEX = (LIMRIGHT-LIMLEFT)/(LIMIT-1)
18. SIZEY = (LIMTOP-LIMBOT)/(LIMIT-1)
19.
20. steps(numSteps,z_r,z_i,c_r,c_i) =
21.     if((numSteps < MAXSTEPS) && (z_r*z_r+z_i*z_i<4),
22.         steps(numSteps+1,
23.         z_r*z_r - z_i*z_i + c_r, 2*z_r*z_i + c_i, c_r, c_i),
24.         numSteps)
25.
26. shape leaves(c_r,c_i) {
27.     FILL[h 214 sat 0.7 b 0.95] # Water
28.     loop i = (LIMIT) [] {
```

International Journal on Integrating Technology in Education (IJITE) Vol.11, No.1, March 2022

```
29.            z_i = (LIMTOP-LIMBOT)*i/(LIMIT-1) + LIMBOT # y
30.         loop j = LIMIT [] {
31.             z_r = (LIMRIGHT-LIMLEFT)*j/(LIMIT-1) + LIMLEFT # x
32.
33.             numSteps = steps(0, z_r, z_i, c_r, c_i)
34.             if(numSteps<MAXSTEPS){
35.                 SQUARE[x z_r
36.                       y z_i
37.                       size SIZEX SIZEY
38.                       h 120 sat 1
39.                       b ((numSteps-1)/(MAXSTEPS-1))]
40.             }
41.         }
42.     }
43. }
```

## 4.2. The crucified

This artwork (figure 7) is based on figure 4, and the author saw some crucified people, and changed seed number, viewport of the complex plane (borders) and some other parameters from listing 4.

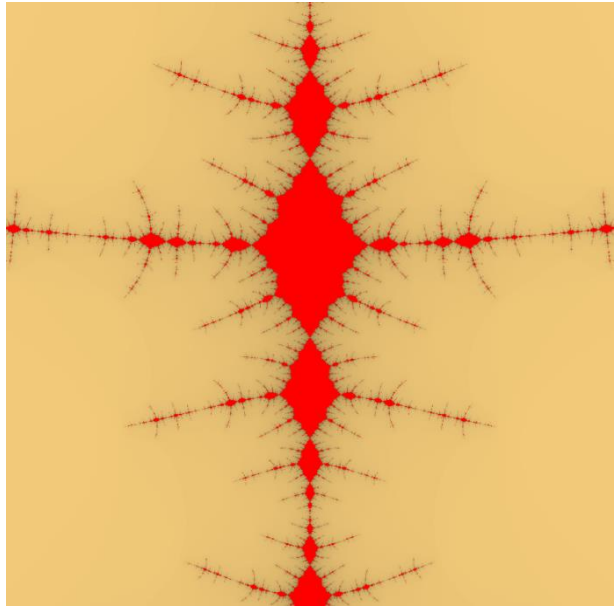

Figure 7. The crucified

Moreover, she had to change the optimization of listing 4 (lines 37–41) since in her new viewport she did not have the symmetry of Figure 4. This was a major problem for her, but she persisted, and with the teacher assistance, she was able to make the necessary modifications. See listing 7.

Listing 7: CFDG source code for The crucified

```
1.  # The crucified
2.  # To generate the image run the following line:
3.  # $ cfdg -b 0 -s 1000 crucified.cfdg crucified.png
4.
5.  startshape crucified(-1.39, 0.0)
6.
7.  LIMIT = 1000 # Image resolution
8.  MAXSTEPS = 200
9.
10. # Borders of the complex plane to show:
```

```
11.  LIMLEFT = -0.02
12.  LIMRIGHT = 0.02
13.  LIMTOP = -0.315
14.  LIMBOT = -0.355
15.
16.  # Width and height of the boxes that will discretize the image:
17.  SIZEX = (LIMRIGHT-LIMLEFT)/(LIMIT-1)
18.  SIZEY = (LIMTOP-LIMBOT)/(LIMIT-1)
19.
20.  steps(numSteps,z_r,z_i,c_r,c_i) =
21.      if((numSteps < MAXSTEPS) && (z_r*z_r+z_i*z_i<4),
22.          steps(numSteps+1,
23.          z_r*z_r - z_i*z_i + c_r, 2*z_r*z_i + c_i, c_r, c_i),
24.          numSteps)
25.
26.  shape crucified(c_r,c_i) {
27.      FILL[b 1 h 0 sat 1] #Red body
28.      loop i = LIMIT [] {
29.          z_i = (LIMTOP-LIMBOT)*i/(LIMIT-1) + LIMBOT # y
30.          loop j = LIMIT [] {
31.              z_r = (LIMRIGHT-LIMLEFT)*j/(LIMIT-1) + LIMLEFT # x
32.
33.              numSteps = steps(0, z_r, z_i, c_r, c_i)
34.              if (numSteps<MAXSTEPS) {
35.                  bright = (1+(1-numSteps)/(MAXSTEPS-1))
36.                  SQUARE[x z_r y z_i size SIZEX SIZEY b bright h 40 sat 0.5]
37.              }
38.          }
39.      }
40.  }
```

## 4.3. Blood sprinkle

The inspiration for this artwork (figure 8) was figure 5, but the author said it was too neat a pattern for a blood sprinkle. So, he thought of superimposing the set of figure 5 three times with different angles.

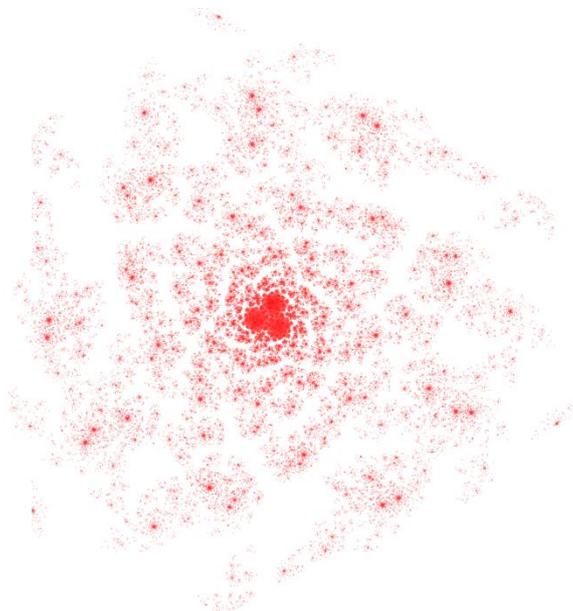

Figure 8.  Blood sprinkle

Triplicating the figure and rotating them was not a big problem (see listing 8, lines 31–37). But the result was not as expected, because in the source code of listing 5, the white dots (squares) are not transparent. So, he had to add a conditional (see listing 8, lines 7 and 44) to prevent the almost white dots (squares) from being generated. The result (figure 8), while equally symmetrical as the original (figure 5), is less neat, as the author targeted.

Listing 8: CFDG source code for Blood sprinkle

```
1.  # Blood sprinkle
2.  # To generate the image run the following line:
3.  # $ cfdg -b 0 -s 1000 blood.cfdg blood.png
4.
5.  LIMIT = 1000 # Image resolution
6.  MAXSTEPS = 150
7.  PROPORTION = 7/10
8.
9.  startshape blood(0.39, -0.252857)
10.
11. SIDE = 0.84 # side of the viewport square
12. CX = 0.21    # center x
13. CY = -0.445714 # center y
14.
15. # Borders of the complex plane to show:
16. LIMLEFT = CX - SIDE/2
17. LIMBOT = CY - SIDE/2
18. LIMRIGHT = CX + SIDE/2
19. LIMTOP = CY + SIDE/2
20.
21. # Width and height of the boxes that will discretize the image:
22. SIZEX = (LIMRIGHT-LIMLEFT)/(LIMIT-1)
23. SIZEY = (LIMTOP-LIMBOT)/(LIMIT-1)
24.
25. steps(numSteps,z_r,z_i,c_r,c_i) =
26.     if((numSteps < MAXSTEPS) && (z_r*z_r+z_i*z_i<4),
27.         steps(numSteps+1,
28.         z_r*z_r - z_i*z_i + c_r, 2*z_r*z_i + c_i, c_r, c_i),
29.         numSteps)
30.
31. shape blood(c_r,c_i){
32.     sprinkle(c_r,c_i)[x -CX -CY]
33.     sprinkle(c_r,c_i)[[r 120 x -CX -CY ]]
34.     sprinkle(c_r,c_i)[[r 240 x -CX -CY ]]
35. }
36.
37. shape sprinkle(c_r,c_i) {
38.     loop i = LIMIT [] {
39.         z_i = (LIMTOP-LIMBOT)*i/(LIMIT-1) + LIMBOT # y
40.         loop j = LIMIT [] {
41.             z_r = (LIMRIGHT-LIMLEFT)*j/(LIMIT-1) + LIMLEFT # x
42.
43.             numSteps = steps(0, z_r, z_i, c_r, c_i)
44.             if (numSteps>PROPORTION*MAXSTEPS) {
45.                 SQUARE[x z_r y z_i size SIZEX SIZEY b 1
46.                     sat ((numSteps-1)/(MAXSTEPS-1))]
47.             }
48.         }
49.     }
50. }
```

## 5. DISCUSSION

After the course, the students expressed being a little surprised by this strange mix-of-maths-and-art sessions, where they were free to experiment and play with the parameters. Moreover, students noted that the fine structures of these images are manifestations of the fact that the smallest variations (mainly the value of the seed) at the beginning of a procedure can result in huge differences later (the different Julia sets are very different from each other), and as [16] tells us, the research of dynamic systems indicates that this is typical of natural processes.

As González Peralta et al. [11] say, there is certainly potential in the inclusion of playful activities in teaching but precautions must be taken to make the sessions useful for the purposes of the curriculum. So, this type of activities should be done mainly in extracurricular spaces, since students have different levels of aptitude and artistic sensitivity.

Artistic activities in general are very enriched thanks to the computer offering the possibility of experimentation, since one can check the influence of parameters on the results, one can check the result of the transformations, the limiting values of interactively applied calculations, etc. [3]. "Modern art studies have shown, however, that meeting the classical definition of beauty is not in itself sufficient to create a work of art. In addition, there must be something to stimulate interest, demand involvement, and motivate further thoughts [3, p. 184]". So, it is not enough to "create" complicated fractal images that are aesthetically beautiful, but they must have a more transcendent meaning.

So, as Sethi and Subramoniam [40] claim, this proposal is meant to accomplish a type of a holistic understanding of Julia Sets and colour variation models, for students to discover meaningful relationships, and develop new knowledge that was difficult to do in the past. Also, this proposal aims to mobilize students in the sense of developing critical capacities with a view to emancipation, like call Ferreira and Lessa [1].

Finally, we can echo McCabe and Reisz's [41] documentary and say that we can extract art from mathematics since "Mathematics is an inherent part of nature".

## AUTHOR

**Eduardo Adam Navas-López** is full time Professor at University of El Salvador. Received the B.Sc. Degree in Computer Science from Central American University "José Simeón Cañas" and M.Sc. Degree in Mathematics Education from University of El Salvador. His research interests are Computer Science Education, Mathematics Education, Computer Graphics, Algorithmic Thinking, and Computational Thinking.

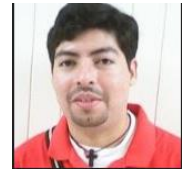